\documentclass[11pt]{amsart}
\usepackage{amsbsy,amssymb,amsmath,amsthm,amscd,amsfonts,latexsym,amstext,delarray,
amsmath,graphicx} 
\usepackage[margin=1in]{geometry}
\usepackage{color}
\input xypic
\usepackage[all]{xy}

\numberwithin{equation}{section}

\def\Q{{\mathbb Q}}
\def\Z{{\mathbb Z}}

\def\C{{\mathbb C}}

\renewcommand{\P}{{\mathbb P}}

\def\bG{{\mathbb G}}
\def\P{{\mathbb P}}

\def\bL{{\mathfb L}}

\def\cA{{\mathcal A}}
\def\cB{{\mathcal B}}
\def\cC{{\mathcal C}}
\def\cD{{\mathcal D}}
\def\cE{{\mathcal E}}

\def\cJ{{\mathcal J}}

\def\cN{{\mathcal N}}
\def\cO{{\mathcal O}}

\def\cV{{\mathcal V}}

\def\cZ{{\mathcal Z}}

\def\bC{{\mathbb C}}

\def\bG{{\mathbb G}}

\def\bL{{\mathbf L}}

\def\bP{{\mathbb P}}
\def\bQ{{\mathbb Q}}

\DeclareMathOperator*{\Hom}{Hom}
\def\Aut{{\rm Aut}}

\def\End{{\rm End}}

\def\Gal{{\rm Gal}}

\def\Hom{{\rm Hom}}
\def\id{{\rm id}}

\def\Spec{{\rm Spec}}

\title{Noncommutative motives and their applications}
\author{Matilde Marcolli and Gon\c{c}alo Tabuada}
\address{Mathematics Department, Caltech, 1200 E. California Blvd. Pasadena, CA 91125, USA}
\email{matilde@caltech.edu}
\urladdr{http://www.its.caltech.edu/~matilde}
\address{Department of Mathematics, MIT, Cambridge, MA 02139, USA}
\email{tabuada@math.mit.edu}
\urladdr{http://math.mit.edu/~tabuada/}
\date{\today}
\thanks{Matilde Marcolli was partially supported by the grants DMS-0901221, DMS-1007207, DMS-1201512, and PHY-1205440. Gon{\c c}alo Tabuada was partially supported by the NEC Award-2742738. \newline 
The authors are very grateful to the organizers Michael Artin, Victor Ginzburg, Catharina Stroppel, Toby Stafford, Michel Van den Bergh, and Efim Zelmanov for kindly giving us the opportunity to present our recent work. They would like also to thank the anonymous referee for comments and corrections.}
\begin{document}

\begin{abstract}
This survey is based on lectures given by the authors during the program ``Noncommutative algebraic geometry and representation theory" at the MSRI in the Spring 2013. It covers the recent work \cite{MaTa2,MaTa1,MaTa3,MaTa4,MaTa6,MaTa5,MaTa7} on noncommutative motives and their applications, and is intended for a broad mathematical audience. In Section 1 we recall the main features of Grothendieck's theory of motives. In Sections 2 and 3 we introduce several categories of noncommutative motives and describe their relation with the classical commutative counterparts. In Section 4 we formulate the noncommutative analogues of Grothendieck's standard conjectures of type C and D, of Voevodsky's smash-nilpotence conjecture, and of Kimura-O'Sullivan finite-dimensionality conjecture. Section 5 is devoted to recollections of the (super-)Tannakian formalism. In Section 6 we introduce the noncommutative motivic Galois (super-)groups and their unconditional versions. 
In Section 7 we explain how the classical theory of (intermediate) Jacobians can be extended to the noncommutative world. Finally, in Section 8 we present some applications to motivic decompositions and to Dubrovin's conjecture.
\end{abstract}
\maketitle
\vskip-\baselineskip
\vskip-\baselineskip

\tableofcontents

\section{Grothendieck's theory of motives}

We recall here the main features of Grothendieck's classical theory of pure
motives, which will be useful when passing to the noncommutative
world. These facts are well-known and we
refer the reader to \cite{Andre,MotConf,Man} for more detailed treatments. Let $k$ be a base field and $F$ a field of coefficients.

\smallskip

Let $\cV(k)$ be the category of smooth projective $k$-schemes. The category of pure motives is obtained from $\cV(k)$ by linearization, idempotent completion, and invertion of the Lefschetz motive.

\subsection{Correspondences}
The linearization of $\cV(k)$ is obtained by replacing the morphisms of schemes with correspondences. Concretely, the correspondences ${\rm Corr}_{\sim,F} (X,Y):=\cZ^{\mathrm{dim}(X)}_{\sim,F}(X \times Y)$ from $X$ to $Y$ are the $F$-linear
combinations of algebraic cycles in $X\times Y$ of codimension equal to $\mathrm{dim}(X)$. This includes the case of ordinary morphisms by
viewing their graphs as correspondences. The composition of correspondences is obtained by pulling back the
cycles to the product $X\times Y\times Z$, taking their intersection
product there, and pushing forward the result to the product $X\times Z$:
\begin{eqnarray}\label{eq:comp}
& {\rm Corr}_{\sim,F}(X,Y)\times {\rm Corr}_{\sim,F}(Y,Z) \rightarrow {\rm Corr}_{\sim,F}(X,Z) & 
(\alpha, \beta)   \mapsto  (\pi_{XZ})_* (\pi_{XY}^*(\alpha) \bullet \pi_{YZ}^*(\beta))\,.
\end{eqnarray}
\subsection{Equivalence relations on algebraic cycles}
One of the important steps in the construction of the category
of pure motives is the choice of an equivalence relation on
algebraic cycles. The usual choices are rational equivalence,  homological equivalence, and numerical equivalence. Rational equivalence depends upon the moving lemma and gives rise to the category of Chow motives. Homological equivalence depends on the choice of a ``good''
cohomology theory (=Weil cohomology theory) and gives rise to the category of homological motives. Numerical equivalence depends only on the intersection product
between algebraic cycles and gives rise to the category of numerical motives. These three equivalence relations are summarized as
follows:
\begin{itemize}
\item A correspondence $\alpha$ from $X$ to $Y$ is rationally trivial, $\alpha \sim_{\mathrm{rat}} 0$, if there exists a $\beta \in \cZ^\ast_{\sim,F}(X\times Y \times \mathbb{P}^1)$ such that $\alpha=\beta(0)-\beta(\infty)$.
\item A correspondence $\alpha$ from $X$ to $Y$ is homologically trivial, $\alpha \sim_{\mathrm{hom}} 0$, if it vanishes under a chosen Weil cohomology theory.
\item A correspondence $\alpha$ from $X$ to $Y$ is numerically trivial, $\alpha \sim_{\mathrm{num}} 0$, if it has a trivial intersection number with every other algebraic cycle.
\end{itemize}
It is well-known that $\sim_{\mathrm{rat}} \neq \sim_{\mathrm{num}}$. The question of whether $\sim_{\mathrm{hom}} = \sim_{\mathrm{num}}$ remains open and is part of an important set of conjectures about motives which will be described below (see \S\ref{sub:standard}).

\smallskip

The category of pure motives has different properties depending on the equivalence relation. 
\subsection{Pure motives}
The symmetric monoidal category of effective pure motives ${\rm Mot}^{\mathrm{eff}}_{\sim, F}(k)$ is defined as follows: the objects are the pairs $(X,p)$ (with $X\in \cV(k)$ and $p$ an idempotent of $\mathrm{Corr}_{\sim,F}(X,X)$), the morphisms are the correspondences
$$ \Hom_{{\rm Mot}^{\mathrm{eff}}_{\sim, F}(k)}((X,p),(Y,q))= p {\rm Corr}_{\sim,F}(X,Y) q\,,$$
the composition is induced from \eqref{eq:comp}, and the symmetric monoidal structure is given by $(X,p)\otimes (Y,q)= (X\times Y, p\times q)$. In what follows we will write $h_F(X)$ instead of $(X,\Delta_X)$.

\smallskip

The effective pure motive $h_F(\P^1)$ decomposes as $h_F^0(\P^1)\oplus h_F^2(\P^1)\simeq 1 \oplus \bL$, where $1=h_F(\Spec(k))$ and $\bL$ is called the Lefschetz motive.

\medskip

The symmetric monoidal category of pure motives ${\rm Mot}_{\sim, F}(k)$ is obtained from  ${\rm Mot}^{\mathrm{eff}}_{\sim, F}(k)$ by formally inverting
the Lefschetz motive. The formal inverse $\bL^{-1}$ is called the Tate motive $\bQ(1)$ (one writes $\Q(n):=\Q(1)^{\otimes n}$). Concretely, in the category of pure motives the objects are triples $(X,p,m):= (X,p)\otimes (\bL^{-1})^{\otimes m}= (X,p)\otimes \bQ(m)$ and the morphisms are given by
$$ \Hom_{{\rm Mot}_{\sim,F}(k)}((X,p,m),(Y,q,n))= p {\rm Corr}^{m-n}_{\sim,F}(X,Y) q\,.$$
The category ${\rm Mot}_{\sim, F}(k)$ is additive. In the case where $m=n$, the direct sum $(X,p,m) \oplus (Y,q,n)$ is defined as $(X \amalg Y, p\oplus q, m)$. The general case reduces to this one using the Lefschetz motive.

\smallskip

Inverting the Lefschetz motive has therefore the effect of introducing arbitrary shifts in the codimension of the algebraic cycles,
instead of using only algebraic cycles of codimension equal to $\mathrm{dim}(X)$. One has a canonical (contravariant) symmetric monoidal functor 
\begin{eqnarray*}
h_F:\cV(k)^{\mathrm{op}} \to {\rm Mot}_{\sim, F}(k) && X \mapsto h_F(X)\,,
\end{eqnarray*}
which sends a morphism $f:X \to Y$ to the transpose of its graph $\Gamma(f)\subset X\times Y$.

\subsection{Chow and homological motives}\label{sub:homological}
The category ${\rm Mot}_{\sim, F}(k)$ with $\sim=\sim_{\mathrm{rat}}$ (resp. $\sim=\sim_{\mathrm{hom}}$) is called the category of Chow motives (resp. homological motives) and is denoted by $\mathsf{Chow}_F(k)$ (resp. by $\mathsf{Hom}_F(k)$).

\subsection{Numerical motives}\label{sub:numerical}
One of the most important results in the theory of pure motives was obtained by
Jannsen \cite{Jan} in the nineties. It asserts that the numerical equivalence relation
is the ``best one" from the point of view of the resulting properties of the category.
More precisely, Jannsen proved that the following conditions are equivalent:
\begin{itemize}
\item ${\rm Mot}_{\sim, F}(k)$ is a semi-simple abelian category;
\item $\mathrm{Corr}_{\sim, F}(X,X)$ is a finite-dimensional
semi-simple $F$-algebra for every $X$;
\item The equivalence relation $\sim$ is equal to $\sim_{\mathrm{num}}$.
\end{itemize}
The category ${\rm Mot}_{\sim, F}(k)$ with $\sim =\sim_{\mathrm{num}}$ is called the category of numerical motives $\mathsf{Num}_F(k)$.

\subsection{Smash-nilpotence}\label{sec:smash}
Voevodsky introduced in \cite{Voe}
the smash-nilpotence equivalence relation on algebraic cycles, $\sim_{\otimes \mathrm{nil}}$, and conjectured the following:

\smallbreak

{\it Conjecture $V(X)$: we have the following equality $\cZ^\ast_{\sim_{\otimes \mathrm{nil}},F}(X) = \cZ^\ast_{\sim_{\mathrm{num}},F}(X)$}. 

\smallbreak

The $\otimes_{\rm nil}$-ideal of an $F$-linear, additive, symmetric monoidal category $\cC$ is defined as 
$$\otimes_{\rm nil}(a,b)=\{g\in {\rm Hom}_{\cC}(a,b)\,|\, g^{\otimes n}=0 \ \mathrm{for} \ n \gg 0 \}\,.$$ 
The quotient functor $\cC \to \cC/\otimes_{\rm nil}$ is $F$-linear, additive, symmetric monoidal,
and conservative. If $\cC$ is idempotent complete, then $\cC/\otimes_{\rm nil}$ is also idempotent complete. One denotes by $\mathsf{Voev}_F(k):= \mathsf{Chow}_F(k)/\otimes_{\rm nil}$ the category of Chow
motives up to smash-nilpotence.

\subsection{All together}
The different categories of pure motives are related by a sequence of $F$-linear, additive, full, symmetric monoidal functors $ \mathsf{Chow}_F(k) \rightarrow  \mathsf{Voev}_F(k) \rightarrow \mathsf{Hom}_F(k) \rightarrow  \mathsf{Num}_F(k)$.

\subsection{Tate motives}
The additive full subcategory of ${\rm Mot}_{\sim, F}(k)$ generated by $\Q(1)$ is called the
category of (pure) Tate motives. This category is independent of the equivalence relation.

\smallskip
\subsection{Weil cohomology theories}

A Weil cohomology theory axiomatizes the properties of
a ``good'' cohomology theory. It consists of a contravariant functor $H^*: \cV(k)^{\mathrm{op}} \to \mathrm{GrVect}(F)$ to $\mathbb{Z}$-graded $F$-vector spaces equipped with the following data:
\begin{itemize}
\item K\"unneth isomorphisms $H^*(X\times Y)\simeq H^*(X)\otimes H^*(Y)$;
\item Trace maps $tr: H^{2\mathrm{dim}(X)}(X)(\mathrm{dim}(X)) \to F$;
\item Cycle maps $\gamma_n : \cZ_{\sim_{\mathrm{rat}},F}^n(X) \to H^{2n}(X)(n)$.
\end{itemize}
One assumes that $\mathrm{dim}\,H^2(\P^1)=1$ and some natural compatibility conditions. Examples of Weil cohomology theories include de Rham,
Betti, \'etale, and cristalline cohomology.

\smallskip

A great deal of difficulty in the theory of pure motives comes from
the poor understanding of the cycle maps. The question of which cohomology classes
are in the range of the $\gamma_n$'s is a notoriously difficult problem
(which includes the Hodge conjecture below). 

\smallskip

The idea of motives can be traced back to Grothendieck's quest for 
a universal cohomology lying behind all the different Weil cohomology theories.

\smallskip
\subsection{Grothendieck's standard conjectures}\label{sub:standard}

Important conjectures relate the
properties of the categories of pure motives with the geometry
of schemes. The standard conjectures are traditionally
labelled as type $C$, $D$, $B$, and $I$. They are summarized
as follows:

\begin{itemize}
\item {\em type $C$ (\text{K\"unneth})}: 
the K\"unneth components of the diagonal $\Delta_X$ are algebraic cycles;
\item {\em type $D$ (\text{Hom=Num})}: the homological and the numerical equivalence relations coincide;
\item {\em type $B$ (\text{Lefschetz})}: the Lefschetz involution $\star_{L,X}$ 
is algebraic (with $\bQ$-coefficients);
\item {\em type $I$ (\text{Hodge})}: the quadratic form defined by the Hodge involution $\star_H$
is positive definite.
\end{itemize}
There are relations between these conjectures: type $B$ and $I$ imply type $D$
and in characteristic zero type $B$ implies all others.  
For our purposes, we will focus on types $C$ and $D$.

\section{From motives to noncommutative motives}

As mentioned above, the origin of pure motives was Grothendieck's quest for a universal cohomology lying behind all the different Weil cohomology theories. The origin of noncommutative
motives is similar. 
In the noncommutative world the basic objects are not schemes but rather dg categories.
Instead of cohomology theories, one has homological type
invariants such as algebraic $K$-theory, cyclic homology (and all its variants), topological Hoschschild homology, etc. In analogy with the commutative world, one can try to identify a
suitable universal invariant lying behind all these different invariants.

\smallskip
\subsection{Dg categories}

A differential graded (=dg) category $\cA$ is a category whose morphism sets 
$\cA(x,y)$ are cochain complexes of $k$-modules ($k$ can more generally be a commutative ring) and whose composition law satisfies the Leibniz rule. A dg functor $F:\cA \to \cB$ is a functor which preserves this extra structure. For further details, we refer the reader to Bondal-Kapranov's pioneering work \cite{BK,BK1} as well as to Keller's ICM survey \cite{ICM}. In what follows, we will denote by $\mathsf{dgcat}(k)$ the category of all (small) 
dg categories and dg functors.
\subsubsection{Perfect complexes}
Dg categories should be understood as ``noncommutative schemes''. The reason for this is that one can canonically associate to every scheme $X$ a dg category, namely the dg category of perfect complexes $\mathrm{perf}_{\mathsf{dg}}(X)$ of $\cO_X$-modules. This dg category enhances the classical derived category of perfect complexes $\mathrm{perf}(X)$ in the sense that the latter is obtained from the former by passing to degree zero cohomology. When $k$ is a field and $X$ is quasi-projective, Lunts and Orlov proved in \cite{LO} that this dg enhancement is in fact ``unique''. This construction gives rise to a well-defined (contravariant) symmetric monoidal functor 
\begin{eqnarray*}
\cV(k)^{\mathrm{op}} \rightarrow \mathsf{dgcat}(k) && X \mapsto \mathrm{perf}_{\mathsf{dg}}(X)\,.
\end{eqnarray*}
An arbitrary dg category should be then considered as the dg category of perfect complexes of an hypothetical ``noncommutative scheme''.
\subsubsection{Saturated dg categories}
Kontsevich introduced in \cite{Kontsevich1,Kontsevich2,Kontsevich3} 
a particular class of dg categories
whose properties closely resemble those of
perfect complexes on smooth proper schemes.
These are called saturated dg categories. Concretely, a dg category $\cA$ is saturated if it is perfect as
a bimodule over itself and if for any two objects $x,y$ we have $\sum_i \mathrm{rank}\, H^i\cA(x,y)< \infty$. A $k$-scheme $X$ is smooth and proper if and only if the associated dg category $\mathrm{perf}_{\mathsf{dg}}(X)$ is saturated.

\smallskip

As mentioned by Kontsevich in \cite{Kontsevich1}, other examples of saturated dg categories arise from representation theory of (finite) quivers and from deformation by quantization.

\subsection{Morita equivalences}

A dg functor $F:\cA \to \cB$ is called a Morita equivalence if the 
restriction of scalars between derived categories $\cD(\cB)\stackrel{\sim}{\to} \cD(\cA)$ 
is an equivalence of (triangulated) categories.

\smallskip

All the classical invariants, like algebraic $K$-theory, cyclic homology (and all its variants), topological Hochschild homology, etc, are Morita invariant in the sense 
that they send Morita equivalences 
to isomorphisms. It is then natural to consider dg categories up to Morita equivalence. 

\smallskip

As proved in \cite{IMRN}, the category $\mathsf{dgcat}(k)$ carries a Quillen model structure whose weak equivalences are the Morita equivalences. Let us denote by $\mathsf{Hmo}(k)$ the homotopy category hence obtained. This category is symmetric monoidal and, as shown in \cite{CisTab}, the saturated dg categories can be characterized as the dualizable (or rigid) objects of $\mathsf{Hmo}(k)$.

\smallskip

Bondal and Van den Bergh proved in \cite{BV} that for every quasi-compact quasi-separated $k$-scheme $X$ the dg category $\mathrm{perf}_{\mathsf{dg}}(X)$ is isomorphic in $\mathsf{Hmo}(k)$ to a dg $k$-algebra with bounded cohomology.

\subsubsection{Bondal-Kapranov's pre-triangulated envelope}
Using Òone-sided twisted complexesÓ, Bondal and Kapranov constructed in \cite{BK} the pre-triangulated envelope of every dg category $\cA$. Intuitively speaking, one formally adds to $\cA$ (de-)suspensions, cones, cones of morphisms between cones, etc. Making use of the Morita model structure, this construction can be conceptually understood as a functorial fibrant resolution functor; consult \cite{IMRN} for details.

\subsubsection{Drinfeld's DG quotient}
A very useful operation on triangulated categories is the Verdier quotient. Via a very elegant construction (reminiscent from Dwyer-Kan localization), Drinfeld \cite{Drinfeld} lifted this operation to the setting of dg categories. Although very elegant, this construction didn't seem to satisfy any obvious universal property. The Morita model structure changed this state of affairs by allowing the characterizing of Drinfeld's DG quotient as an homotopy cofiber construction; consult \cite{Drinfeld-Tab} for details.

\smallskip
\subsection{Additive invariants}\label{sub:additive}
Given a dg category $\cA$, let $T(\cA)$ be the dg category of pairs $(i,x)$, where $i\in \{1,2\}$ and $x \in \cA$. The complex of morphisms in $T(\cA)$ from $(i,x)$ to $(i',x')$ is given by $\cA(x,x')$ if $i \leq i'$ and is zero otherwise. Composition is induced by $\cA$. Intuitively speaking, $T(\cA)$ ``dg categorifies'' the notion of upper triangular matrix.  Note that we have two inclusion dg functors $i_1: \cA \hookrightarrow T(\cA)$ and $i_2: \cA \hookrightarrow T(\cA)$. Let $E:\mathsf{dgcat}(k)\to \mathsf{A}$ be a functor with values in an additive category. We say that $E$ is an additive invariant if it satisfies the following two conditions:
\begin{itemize}
\item it sends Morita equivalences to isomorphisms;
\item given any dg category $\cA$, the inclusion dg functors induce an isomorphism
$$ [E(i_1)\,\,E(i_2)]:E(\cA) \oplus E(\cA) \stackrel{\sim}{\longrightarrow} E(T(\cA))\,.$$
\end{itemize}
Thanks to the work of Blumberg-Mandell, Keller, Quillen, Schlichting, Thomason-Trobaugh, Waldhausen, and others (see \cite{BluMa,Kell, Quillen, Schl, MacLane, ThoTro, Wald}), all the above invariants are additive.

\smallskip

The universal additive invariant was constructed in \cite{IMRN} as follows: consider the additive symmetric monoidal category $\mathsf{Hmo}_0(k)$
whose objects are the dg categories and whose morphisms are given by $\Hom_{\mathsf{Hmo}_0(k)}(\cA,\cB) := K_0 \mathrm{rep}(\cA,\cB)$, where ${\rm rep}(\cA,\cB) \subset \cD(\cA^{\mathrm{op}}\otimes^{\bL} \cB)$ stands for the full triangulated subcategory of those $\cA\text{-}\cB$ bimodules which are perfect as right $\cB$-modules. The composition law and the symmetric monoidal structure are induced from $\mathsf{Hmo}(k)$. As explained in loc. cit., the canonical composed symmetric monoidal functor 
\begin{equation}\label{eq:universal-U}
U:\mathsf{dgcat}(k) \longrightarrow \mathsf{Hmo}(k) \longrightarrow \mathsf{Hmo}_0(k)
\end{equation}
is the universal additive invariant, i.e. pre-composition with $U$ induces a bijection between additive functors $\overline{E}: \mathsf{Hmo}_0(k) \to \mathsf{A}$ and additive invariants $E:\mathsf{dgcat}(k) \rightarrow \mathsf{A}$. This suggests that $\mathsf{Hmo}_0(k)$ is the correct place where noncommutative motives should reside. Let us denote by $\mathsf{Hmo}_0(k)_F$ its $F$-linearization ($F$ can more generally be a commutative ring).

\subsection{Computations}\label{sub:computations}
In order to gain some sensibility with \eqref{eq:universal-U}, we recall some computations:
\begin{itemize}
\item Let $X$ be a smooth projective $k$-scheme whose derived category $\mathrm{perf}(X)$ admits a full exceptional collection of length $n$ (see \S\ref{sub:exceptional}). In this case, $U(\mathrm{perf}_{\mathsf{dg}}(X)) $ identifies with the direct sum of $n$ copies of $U(k)$; see~\cite{MaTa4}.
\item Let $X$ be a quasi-compact quasi-separated $k$-scheme, $A$ a sheaf of Azumaya algebras over $X$ of rank $r$, and $\mathrm{perf}_{\mathsf{dg}}(A)$ the associated dg category. When $1/r$ belongs to the commutative ring $F$, $U(\mathrm{perf}_{\mathsf{dg}}(A))_F$ identifies with $U(\mathrm{perf}_{\mathsf{dg}}(X))_F$; see \cite{TV}.
\item Let $A$ be a central simple $k$-algebra of degree $d:=\sqrt{\mathrm{dim}(A)}$ and $SB(A)$ the associated Severi-Brauer variety. In this case, $U(\mathrm{perf}_{\mathsf{dg}}(SB(A)))$ identifies with the following direct sum $U(k) \oplus U(A) \oplus U(A)^{\otimes 2} \oplus \cdots \oplus U(A)^{\otimes (d-1)}$; see \cite{MarTab3}.
\item Let $k$ be a perfect field, $A$ a finite-dimensional $k$-algebra of finite global dimension, and $J(A)$ its Jacobson radical. In this case, we have $U(A) \simeq U(A/J(A))$; see \cite{TV}.
\item Let $A$ and $B$ be two central simple $k$-algebras and $[A]$ and $[B]$ their Brauer classes. In this case, $U(A)$ and $U(B)$ are isomorphic if and only if $[A]=[B]$; see \cite{Twisted}.
\end{itemize}

\subsection{Noncommutative Chow motives}

Kontsevich introduced in \cite{Kontsevich1} the symmetric monoidal category of  noncommutative Chow motives $\mathsf{NChow}(k)_F$. It can be described as the idempotent completion $(-)^\natural$ of the full subcategory of $\mathsf{Hmo}_0(k)_F$ given by the saturated dg categories. Concretely, the objects are of the pairs $(\cA,e)$ (with $\cA$ a saturated dg category and $e$ an idempotent), the morphisms are given by the noncommutative correspondences
$$\Hom_{\mathsf{NChow}(k)_F}((\cA,e),(\cB,e')) :=e K_0 \mathrm{rep}(\cA,\cB) e' \simeq e K_0(\cA^{\mathrm{op}} \otimes^\bL \cB)_F e'\,,$$ 
the composition law is induced by the derived tensor product of bimodules, and the symmetric monoidal structure is induced by the derived tensor product of dg categories.

\subsubsection{Fundamental theorems}
The fundamental theorems in algebraic K-theory and periodic cyclic homology, proved respectively by Weibel \cite{Weibel} and Kassel \cite{Kassel}, are of major importance. Their proofs are not only very different but also quite involved. The category $\mathsf{NChow}(k)_F$ allowed a simple, unified and conceptual proof of these fundamental theorems; consult \cite{Fundamental} for details.

\smallskip
\subsection{A bridge from Chow to noncommutative Chow motives}\label{sub:bridge}

Noncommutative motives should, in a suitable sense, contain
the category of motives. This idea was made precise in \cite{Tab2} (following Kontsevich's original insight \cite{Kontsevich3}). The precise
statement is the existence of a $\bQ$-linear additive fully-faithful symmetric monoidal functor $R$ making the following diagram commute
\begin{equation}\label{eq:diag-main}
\xymatrix{
\cV(k)^{\mathrm{op}} \ar[d]_-{h_\bQ} \ar[rr]^-{X \mapsto \mathrm{perf}_{\mathsf{dg}}(X)} && \mathsf{dgcat}(k) \ar[d]^-U \\
\mathsf{Chow}_\bQ(k) \ar[d]_-\pi && \mathsf{Hmo}_0(k) \ar[d]^-{(-)_\bQ} \\
\mathsf{Chow}_\bQ(k)/_{\!\!-\otimes \bQ(1)} \ar[rr]_-R && \mathsf{NChow}(k)_\bQ \subset \mathsf{Hmo}_0(k)_\bQ^\natural\,,
}
\end{equation}
where $\mathsf{Chow}_\bQ(k)/_{\!\!-\otimes \bQ(1)}$ stands for the orbit category. This bridge opens new horizonts and opportunities of research by enabling the interchange of results between the commutative and the noncommutative world. This yoga was developed in \cite{Tab2} in what regards Schur and Kimura-O'Sullivan finite-dimensionality (see \S\ref{sub:finite} below), motivic measures, and motivic zeta functions.

\smallskip

The above diagram \eqref{eq:diag-main} holds more generally with $\bQ$ replaced by any field $F$ of characteristic zero.

\subsubsection{Orbit categories} 
Given an $F$-linear, additive, symmetric monoidal category $\cC$ and a $\otimes$-invertible object $\cO \in \cC$, the orbit category $\cC/_{\!\!-\otimes \cO}$ has the same objects as $\cC$ and morphisms
$$\Hom_{\cC/_{\!\!-\otimes \cO}} (a,b):= \bigoplus_{i\in \Z} \Hom_{\cC} (a, b\otimes \cO^{\otimes i})\,.$$ 
The composition law and the symmetric monoidal structure are induced from $\cC$. By construction, the orbit category comes equipped with a canonical symmetric monoidal projection functor $\pi:\cC \to \cC/_{\!\!-\otimes \cO}$. Moreover, $\pi$ is endowed with a $2$-isomorphism $\pi \circ (- \otimes \cO) \stackrel{\sim}{\Rightarrow} \pi$ and is $2$-universal among all such functors.

\section{Categories of noncommutative motives}

\subsection{Periodic cyclic homology as ``noncommutative de Rham cohomology''}\label{sub:HP}
Connes' periodic cyclic homology extends naturally from $k$-algebras to dg categories giving thus rise to a functor $HP^\pm:\mathsf{dgcat}(k) \rightarrow \mathrm{sVect}(k)$ to super $k$-vector spaces. In the case of a smooth $k$-scheme $X$ (with $k$ of characteristic zero), the Hochschild-Konstant-Rosenberg
theorem show us that 
\begin{equation}\label{eq:HKR}
HP^\pm(\mathrm{perf}_{\mathsf{dg}}(X)) \simeq HP^\pm(X) \simeq (\bigoplus_{n \,\mathrm{even}} H^n_{dR}(X), 
        \bigoplus_{n \, \mathrm{odd}} H^n_{dR}(X))\,.
\end{equation}
For this reason $HP^\pm$ is considered the noncommutative analogue of de Rham cohomology. For further details on this viewpoint, we invite the reader to consult Kaledin's ICM address \cite{Kaledin}.

\smallskip

As proved in \cite[Theorem~7.2]{MaTa3}, $HP^\pm$ induces symmetric monoidal functors
\begin{eqnarray}\label{eq:HP1}
\overline{HP^\pm}: \mathsf{NChow}(k)_F \to {\rm sVect}(F) && \overline{HP^\pm}: \mathsf{NChow}(k)_F \to {\rm sVect}(k)
\end{eqnarray}
under the assumption that $F$ is a field extension of $k$ (left-hand-side) or the assumption that $k$ is a field extension of $F$ (right-hand-side). 

\subsection{Noncommutative homological motives}

Making use of the above ``noncommutative de Rham cohomology'', one defines the symmetric monoidal category of  noncommutative homological motives $\mathsf{NHom}(k)_F$ as the idempotent completion of the quotient category $\mathsf{NChow}(k)_F/{\rm Ker}(\overline{HP^\pm})$. 

\subsection{Noncommutative numerical motives}\label{sub:NCnumerical}

In order to define a category of noncommutative numerical motives one needs to extend to the noncommutative world 
the notion of intersection number. This can be done as follows. Let $(\cA,e)$ and $(\cB,e')$ be two noncommutative Chow motives. Given a noncommutative correspondence $\underline{\mathsf{B}} =e [\sum_i a_i \mathsf{B}_i] e'$ from $(\cA,e)$ to $(\cB,e')$ (recall that the $\mathsf{B}_i$'s are $\cA\text{-}\cB$-bimodules), and a noncommutative correspondence $\underline{\mathsf{B}'}= e' [\sum_j b_j \mathsf{B}'_j] e$ from $(\cB,e')$ to $(\cA,e)$, one defines their intersection number as the following sum
$$ \langle \underline{\mathsf{B}},\underline{\mathsf{B}'} \rangle := \sum_{i,j, n} (-1)^n a_i \cdot b_j \cdot \mathrm{rank}\, HH_n(\cA; \mathsf{B}_i\otimes_{\cB}^{\bL} \mathsf{B}'_j) \,,$$
where $HH_n(\cA; \mathsf{B}_i\otimes_{\cB}^{\bL} \mathsf{B}'_j)$ stands for the $n^{\mathrm{th}}$ Hochschild homology group of $\cA$ with coefficients in the $\cA\text{-}\cA$ bimodule 
$\mathsf{B}_i\otimes_{\cB}^{\bL} \mathsf{B}'_j$. The numerical equivalence relation on noncommutative Chow motives is obtained by declaring a noncommutative correspondence  $\underline{\mathsf{B}}$ to be numerically trivial if
$\langle \underline{\mathsf{B}},\underline{\mathsf{B}'} \rangle = 0$
for all $\underline{\mathsf{B'}}$.
This equivalence relation gives rise to the largest $\otimes$-ideal $\cN$ strictly contained in $\mathsf{NChow}(k)_F$.
The symmetric monoidal category of noncommutative numerical motives $\mathsf{NNum}(k)_F$
is then defined as the idempotent completion of the quotient category $\mathsf{NChow}(k)_F/ \cN$.

\smallskip

As proved in \cite[Theorem~1.12]{MaTa1}, the functor $R$ of diagram \eqref{eq:diag-main} descends to a $\bQ$-linear additive fully-faithful symmetric monoidal functor $R_{\mathsf{N}}:\mathsf{Num}_\Q(k)/_{\!\!-\otimes \Q(1)}\to\mathsf{NNum}(k)_\Q$.

\smallskip

A different numerical equivalence relation on noncommutative Chow
motives (based on the bilinear pairing $\sum_i (-1)^i \mathrm{dim}\, \mathrm{Ext}^i(-,-)$), was proposed by Kontsevich in \cite{Kontsevich1}. As proved in \cite[Theorem~1.1]{MaTa2}, Kontsevich's notion is equivalent to the above one.

\subsection{Semi-simplicity}
As proved in \cite[Theorem~1.1]{MaTa1}\cite[Theorem~4.6]{MaTa3}, Jannsen's result (see \S\ref{sub:numerical}) holds also in the noncommutative world. Concretely, under the assumption that $k$ and $F$ have same characteristic, $\mathsf{NNum}(k)_F$ is abelian semi-simple. This was conjectured by Kontsevich in \cite{Kontsevich1}. Moreover, if $\cJ$ is a $\otimes$-ideal in $\mathsf{NChow}(k)_F$ for which the idempotent completion of $\mathsf{NChow}(k)_F/\cJ$ is abelian semi-simple, then $\cJ$ agrees with $\cN$.

\smallskip

As explained in \cite[Corollary~1.1]{MaTa1}, the semi-simplicity of $\mathsf{NNum}(k)_F$ (with $k$ of characteristic zero) combined with the functor $R$ of diagram \eqref{eq:diag-main} gives rise to an alternative proof of Jannsen's result. 

\subsection{Smash-nilpotence in the noncommutative world} 

Recall from \S\ref{sec:smash} the definition of the $\otimes_{\mathrm{nil}}$-ideal. One denotes by $\mathsf{NVoev}(k)_F:=\mathsf{NChow}(k)_F/\otimes_{\rm nil}$ the category of noncommutative Chow motives up to smash-nilpotence. As proved in \cite[Proposition~3.1]{MaTa7}, the functor $R$ of diagram \eqref{eq:diag-main} descends also to a symmetric monoidal fully-faithful functor $R_{\mathsf{V}}:\mathsf{Voev}_\Q(k)/_{\!\!-\otimes \Q(1)}\to\mathsf{NVoev}(k)_\Q$.

\smallskip

\subsection{All together} 
The categories of noncommutative motives are related by $F$-linear, additive, full, symmetric monoidal functors $ \mathsf{NChow}(k)_F \rightarrow \mathsf{NVoev}(k)_F \rightarrow \mathsf{NHom}(k)_F \rightarrow \mathsf{NNum}(k)_F$. 

\smallskip

Given a saturated dg category $\cA$, we will denote by $\sim_{\otimes \mathrm{nil}}, \sim_{\mathrm{hom}}, \sim_{\mathrm{num}}$ the equivalence relations on $\mathrm{Hom}_{\mathsf{NChow}(k)_F}(U(k)_F,U(\cA)_F) \simeq K_0(\cA)_F$ induced by the above functors.

\subsection{Noncommutative Artin motives}
The category of Artin motives $\mathsf{AM}_F(k)$ is by definition the smallest additive rigid idempotent complete full subcategory of $\mathsf{Chow}_F(k)$ containing the finite {\'e}tale $k$-schemes. One defines the category of noncommutative Artin motives $\mathsf{NAM}(k)_\bQ$ as the image of $\mathsf{AM}_\bQ(k)$ under the functors $R\circ \pi$ of diagram \eqref{eq:diag-main}. As proved in \cite[Theorem 1.7]{MaTa6}, this latter category is independent of the equivalence relation. 

\section{Conjectures in the noncommutative world}\label{sec:conjectures}
Let $\cA$ be a saturated dg category.
\subsection{Standard conjectures}\label{StConjNCsec}
The authors introduced in \cite{MaTa3} the noncommutative analogues of Grothendieck's standard conjectures of type $C$ and $D$. 

\smallskip

{\it Conjecture $C_{NC}(\cA)$:  the K\"unneth projectors 
\begin{eqnarray*}
\pi_\cA^+: HP^\pm(\cA) \to  HP^+(\cA)
\to  HP^\pm(\cA)&&\pi_\cA^-: HP^\pm (\cA) \to  HP^-(\cA)
\to  HP^\pm(\cA)
\end{eqnarray*}
are algebraic, i.e. $ \pi_{\cA}^\pm= \overline{HP^\pm}(\underline{\pi}_{\cA}^\pm)$ for noncommutative correspondences $\underline{\pi}_{\cA}^\pm$.} 

\smallskip

A weaker version of the standard conjecture of type $C$ (K\"unneth) is
the sign conjecture $C^+(X)$: {\it the K\"unneth projectors $ \pi_X^+=\sum_i \pi^{2i}_X$ and $ \pi_X^-=\sum_i \pi^{2i+1}_X$  are algebraic}. The restriction of $C_{NC}$ to the commutative world is weaker than the sign conjecture in the sense that 
$C^+(X) \Rightarrow C_{NC}(\mathrm{perf}_{\mathsf{dg}}(X))$.

\smallskip

{\em Conjecture $D_{NC}(\cA)$: one has the following equality $K_0(\cA)_F/\!\!\sim_{\mathrm{hom}} = K_0(\cA)_F/\!\!\sim_{\mathrm{num}}$}. 

\smallskip

Similarly, the restriction of $D_{NC}$ to the commutative world is weaker than the standard conjecture of type $D$ (Hom=Num) in the sense that  $D(X) \Rightarrow D_{NC}(\mathrm{perf}_{\mathsf{dg}}(X))$.

\subsection{Smash-nilpotence conjecture}
Voevodsky's nilpotence conjecture (see \S\ref{sec:smash}) was extended in \cite{MaTa7} to the noncommutative world as follows:

\smallskip

{\it Conjecture $V_{NC}(\cA)$: we have the following equality $K_0(\cA)_F/\!\!\sim_{\otimes_{\mathrm{nil}}} = K_0(\cA)_F /\!\!\sim_{\mathrm{num}}$}.

\smallskip

As proved in \cite[Theorem~4.1]{MaTa7}, the restriction of $V_{NC}$ to the commutative world is equivalent to Voevodsky's smash-nilpotence conjecture in the sense that $V(X) \Leftrightarrow V_{NC}(\mathrm{perf}_{\mathsf{dg}}(X))$. This suggests that instead of attacking Voevodsky's conjecture $V$, one should alternatively attack conjecture $V_{NC}$ (using noncommutative tools). The authors are currently developing this approach.

\subsection{Finite-dimensionality conjecture}\label{sub:finite}
Let $F$ be a field of characteristic zero and $\cC$ an $F$-linear, idempotent complete, symmetric monoidal category. An object $a \in \cC$ is called even (resp. odd) dimensional if $\wedge^n(a) =0$ (resp. $\mathrm{Sym}^n(a) =0$) for some $n>0$. An object $a \in \cC$ is called finite-dimensional if $a=a_+ \oplus a_-$, with $a_+$ (resp. $a_-$) even (resp. odd) dimensional. Kimura \cite{Kim} and  O'Sullivan \cite{OSul} conjectured the following:

\smallskip

{\it Conjecture $KS(X)$: the Chow motive $h_F(X)$ is finite-dimensional.}

\smallskip

This conjecture was extended in \cite{MaTa7} to the noncommutative world as follows:

\smallskip

{\it Conjecture $KS_{NC}(\cA)$: the noncommutative Chow motive $U(\cA)_F$ is finite-dimensional.}

\smallskip

The restriction of $KS_{NC}$ to the noncommutative world is weaker than the Kimura-O'Sullivan finite-dimensionality conjecture in the sense that $KS(X) \Rightarrow KS_{NC}(\mathrm{perf}_{\mathsf{dg}}(X))$.
\smallskip

Under the assumption that $k$ is a field extension of $F$ (or vice-versa), it was proved in \cite[Theorem 4.1]{MaTa7} that conjectures $V_{NC}((\cA^{\mathrm{op}})^{\otimes n} \otimes \cA^{\otimes n}), n \geq 1$, combined with conjecture $C_{NC}(\cA)$, imply conjecture $KS_{NC}(\cA)$. Moreover, if conjecture $KS_{NC}$ holds for every saturated dg category and all symmetric monoidal functors $\mathsf{NChow}(k)_F \to \mathrm{sVect}(K)$ (with $K$ a field extension of $F$) factor through $\mathsf{NNum}(k)_F$, then conjecture $V_{NC}$ also holds for every saturated dg category. 

\section{(Super-)Tannakian formalism}

\subsection{Tannakian categories} Let $(\cC,\otimes, {\bf 1})$ be an $F$-linear, abelian, symmetric monoidal category. In particular, we have commutativity and $\otimes$-unit constraints
\begin{equation*}
 \tau_{a,b}: a\otimes b \stackrel{\sim}{\rightarrow} b \otimes a  \ \ \ \ \ \ \ell_a: a \stackrel{\sim}{\to} a \otimes {\bf 1} \ \ \ \ \ \ 
r_a: {\bf 1}\otimes a \stackrel{\sim}{\rightarrow} a
\end{equation*}
and the following equality holds $\tau_{b,a}\circ  \tau_{a,b} = \id_{a\otimes b}$.
The category $\cC$ is called rigid if there exists a duality functor
$(-)^\vee: \cC \to \cC^{\mathrm{op}}$, evaluation maps $\epsilon: a\otimes a^\vee \to {\bf 1}$, and coevaluation maps
$\eta: {\bf 1} \to a^\vee \otimes a$, for which the following composition is equal to the identity
$$ a \stackrel{\ell_a}{\longrightarrow} a \otimes {\bf 1} \stackrel{\id_a \otimes \eta}{\longrightarrow} a\otimes a^\vee \otimes a \stackrel{\epsilon \otimes \id_a}{\longrightarrow} {\bf 1} \otimes a \stackrel{r_a}{\longrightarrow} a\,.$$ 
The categorical trace of an endomorphism $g:a \to a$ is defined as $tr(g) = \epsilon \circ \tau_{a^\vee \otimes a} \circ (\id_{a^\vee} \otimes g)  \circ \eta$. The number $\dim(a):= tr(\id_a)$ is called the dimension or Euler characteristic of $a$.

\smallskip

A category $\cC$ with the above properties, and with $\End({\bf 1})=F$, is called 
Tannakian if there exists an exact faithful symmetric monoidal functor $\omega: \cC \to \mathrm{Vect}(K)$ with values in a category of $K$-vector spaces (with $K$ a field extension of $F$). The functor $\omega$ is called a fiber functor. When this holds with $K=F$, $\cC$ is called a  neutral Tannakian category.

\smallskip

If $\cC$ is a neutral Tannakian category, then there is a $\otimes$-equivalence of
categories $\cC \simeq {\rm Rep}_F(\mathrm{Gal}(\cC))$. The affine group scheme 
$\Gal(\cC)$ is given by the $\otimes$-automorphisms $\underline{\mathrm{Aut}}^\otimes(\omega)$
of the fiber functor $\omega$.

\smallskip

Tannakian categories in characteristic zero were characterized by Deligne \cite{Del1} as follows: an $F$-linear, abelian, rigid symmetric monoidal category $\cC$, with $\End({\bf 1})=F$, is
Tannakian if and only if $\mathrm{dim}(a) \geq 0$ for 
all objects $a \in \cC$.

\subsection{Super-Tannakian categories}

An $F$-linear, abelian, rigid symmetric monoidal category $\cC$, 
with ${\rm End}({\bf 1})=F$, is called super-Tannakian if there exists a super-fiber functor $\omega: \cC \to {\rm sVect}(K)$ with values in a category of super $K$-vector spaces (with $K$ a field extension of $F$). When this holds with $K=F$, $\cC$ is called a neutral super-Tannakian category.

\smallskip

If $\cC$ is a neutral super-Tannakian category, $\omega$ induces a $\otimes$-equivalence  between $\cC$ and the category ${\rm Rep}_F({\rm sGal}(\cC),\epsilon)$ of super-representations 
of the affine super-group scheme ${\rm sGal}(\cC) := \underline{\Aut}^\otimes(\omega)$ (the super-structure is given by the parity automorphism
$\epsilon$).

\smallskip

Super-Tannakian categories were also characterized by Deligne \cite{Del2} as follows: an $F$-linear, abelian, rigid symmetric monoidal category $\cC$, 
with ${\rm End}({\bf 1})=F$, is
super-Tannakian if and only if it is Schur-finite. When $F$ is algebraically closed, $\cC$ is
neutral super-Tannakian if and only if it is Schur-finite.

\subsubsection{Schur-finiteness}
Let $\cC$ be a category as above, $S_n$ the symmetric
group on $n$ symbols, and $\Q[S_n]$ the associated group ring. Every partition $\lambda$ of $n$ gives rise to an idempotent $e_\lambda \in \Q[S_n]$ and to a Schur functor $S_\lambda: \cC \to \cC, a \mapsto e_\lambda(a^{\otimes n})$. The category $\cC$ is called Schur-finite if all its objects are annihilated by some Schur functor.

\smallskip
\subsection{Motivic Galois groups}\label{sub:motivic}
Deligne's characterization of Tannakian categories is not satisfied in the case of $\mathsf{Num}_F(k)$ because
$\mathrm{dim}(h_F(X))$ is equal to the Euler characteristic $\chi(X)$ of the $k$-scheme $X$ which can be negative. Jannsen proved in \cite{Jan} that if the standard conjecture of 
type $C$ holds, then one can modify the commutativity constraints $\tau_{X,Y}$
using the algebraic cycles coming from the K\"unneth components
of the diagonal. This has the effect of correcting the negative signs of the Euler characteristic. Let $\mathsf{Num}_F^\dagger(k)$ be the Tannakian category hence obtained. If the standard conjecture of type $D$ also holds, then $\mathsf{Num}_F^\dagger(k)$ is a neutral Tannakian category and every Weil cohomology theory $H^*$ is a fiber functor. Under these assumptions, one obtains a  group scheme $\mathrm{Gal}(\mathsf{Num}^\dagger_F(k))$ called the motivic Galois group.

\subsection{Motivic Galois super-groups}\label{sub:supermotivic}

In contrast with \S\ref{sub:motivic}, Deligne's characterization of super-Tannakian categories is satisfied in the case of $\mathsf{Num}_F(k)$. When $F$ is algebraically closed, $\mathsf{Num}_F(k)$ is then a neutral Tannakian category. As a consequence, one obtains a super-group scheme $\mathrm{sGal}(\mathsf{Num}_F(k))$ called the motivic Galois super-group.

\section{Noncommutative motivic Galois (super-)groups} 

Let $k$ be a field of characteristic zero and $F$ a field extension of $k$.

\smallskip

Assuming conjectures $C_{NC}$  and $D_{NC}$, it was proved in \cite[Theorem 1.6]{MaTa3} that the category $\mathsf{NNum}^\dagger(k)_F$ (obtained from $\mathsf{NNum}(k)_F$ by modifying the commutativity constrains) is neutral Tannakian. An explicit fiber functor is given by periodic cyclic homology. The associated group scheme $\mathrm{Gal}(\mathsf{NNum}^\dagger(k)_F)$ is called the noncommutative motivic Galois group.

\smallskip

As proved in \cite[Theorem~1.2]{MaTa3}, the category 
$\mathsf{NNum}(k)_F$ is super-Tannakian. When $F$ is algebraically closed, $\mathsf{NNum}(k)_F$ is neutral super-Tannakian. Under these assumptions, one obtains a super-group scheme $\mathrm{sGal}(\mathsf{NNum}(k)_F)$ called the noncommutative motivic Galois super-group.

\subsection{Comparison morphisms}
Assuming conjectures $C,D, C_{NC},D_{NC}$, we have well-defined (noncommutative) motivic Galois (super-)groups. As proved in \cite[Theorem 1.7]{MaTa3}, the composed functor
\begin{equation}\label{eq:composed}
\mathsf{Num}_k(k) \stackrel{\pi}{\longrightarrow} \mathsf{Num}_k(k)/_{\!\!-\otimes \bQ(1)} \stackrel{R_{\mathsf{N}}}{\longrightarrow} \mathsf{NNum}(k)_k
\end{equation}
gives rise to faithfully-flat comparison morphisms
\begin{eqnarray}
\mathrm{Gal}(\mathsf{NNum}^\dagger(k)_k) & \relbar\joinrel\twoheadrightarrow&  {\rm Ker}(t: {\rm Gal}(\mathsf{Num}^\dagger_k(k))
\twoheadrightarrow \bG_m) \label{eq:morph1} \\
\mathrm{sGal}(\mathsf{NNum}(k)_k) &\relbar\joinrel\twoheadrightarrow & {\rm Ker}(t: {\rm sGal}(\mathsf{Num}_k(k))
\twoheadrightarrow \bG_m)\,,
\label{eq:morph2}
\end{eqnarray}
where $\bG_m$ is the multiplicative group scheme and $t$ is induced by the inclusion of  Tate motives in $\mathsf{Num}_k(k)$. These comparison morphisms were suggested by Kontsevich in \cite{Kontsevich1}. Intuitively speaking, they show us that the $\otimes$-symmetries of the commutative world which can be lifted to the noncommutative world are precisely those that become trivial when restricted to Tate motives.

\smallskip

The proof of \eqref{eq:morph1}-\eqref{eq:morph2} makes use of Deligne-Milne's theory of Tate triples (see \cite{DelMil}), of a suitable extension of this theory to the super-Tannakian setting, and of Milne's work \cite{Mil} on quotients of Tannakian categories. The key step is the description of the right-hand-side of \eqref{eq:morph1} (resp. of \eqref{eq:morph2}) as the Galois group (resp. super-group) of the orbit category of $\mathsf{Num}^\dagger_k(k)$ (resp. of $\mathsf{Num}_k(k)$).

\smallskip

It is unclear at the moment if the kernel of these comparison morphisms is non-trivial. This problem is related to the existence of ``truly noncommutative numerical motives'', i.e. objects of $\mathsf{NNum}(k)_k$ that are not in the essential image of \eqref{eq:composed}.

\subsection{Unconditional version}
The functors \eqref{eq:HP1} descend to symmetric monoidal functors
\begin{equation}\label{eq:HP2}
\overline{HP^\pm}: \mathsf{NHow}(k)_F \longrightarrow {\rm sVect}(K)\,.
\end{equation}
Here, $K=F$ when $F$ is a field extension of $k$ and $K=k$ when $k$ is a field extension of $F$. Let $\mathsf{NHom}(k)_F^\pm$ be the full subcategory of $\mathsf{NHom}(k)_F$ consisting of those noncommutative homological motives for which the K\"unneth projectors are algebraic. By changing the commutativity constraints of this latter category, one obtains a rigid symmetric monoidal category $\mathsf{NHom}^\dagger(k)_F^\pm$ and an $F$-linear symmetric monoidal functor $\mathsf{NHom}^\dagger(k)_F^\pm \to \mathrm{Vect}(K)$. Making use of Andr{\'e}-Kahn's techniques \cite{AnKa,AnKa2}, the authors showed in \cite[\S1]{MaTa7} that the associated category $\mathsf{NNum}^\dagger(k)_F^\pm$ is Tannakian, semi-simple, and that the canonical functor $\mathsf{NHom}^\dagger(k)_F^\pm \to 
\mathsf{NNum}^\dagger(k)_F^\pm$ admits a $\otimes$-section $s^{NC}$ (unique up to
conjugation by a $\otimes$-isomorphism). One obtains in this way a fiber functor
$$ \omega: \mathsf{NNum}^\dagger(k)_F^\pm \stackrel{s^{NC}}{\longrightarrow} 
 \mathsf{NHom}^\dagger(k)_F^\pm \longrightarrow {\rm Vect}(K). $$
The associated group scheme ${\rm Gal}(\mathsf{NNum}^\dagger(k)_F^\pm)$ is called the unconditional noncommutative motivic Galois group. As proved in \cite[Theorem~1.7]{MaTa7}, we have a faithfully-flat comparison morphism
\begin{eqnarray}\label{eq:uncond}
\mathrm{Gal}(\mathsf{NNum}^\dagger(k)_k^\pm) & \relbar\joinrel\twoheadrightarrow&  {\rm Ker}(t: {\rm Gal}(\mathsf{Num}^\dagger_k(k)^\pm)
\twoheadrightarrow \bG_m )\,.
\end{eqnarray}
Assuming conjectures $C,D$ and $C_{NC},D_{NC}$, the unconditional noncommutative motivic Galois group agree with the conditional one $\mathrm{Gal}(\mathsf{NNum}^\dagger(k)_k)$ and \eqref{eq:uncond} identifies with \eqref{eq:morph2}.

\smallskip
\subsection{Base-change}
As proved in \cite[Theorem 1.9]{MaTa6}, one has the following short exact sequence
\begin{equation}\label{eq:base-change1}
1 \longrightarrow {\rm Gal}(\mathsf{NNum}^\dagger(\bar k)_F)\stackrel{I}{\longrightarrow}
{\rm Gal}(\mathsf{NNum}^\dagger(k)_F)\stackrel{P}{\longrightarrow} {\rm Gal}(\bar k/k) \longrightarrow 1\,.
\end{equation}
Some explanations are in order: $I$ is induced by base-change $-\otimes_k \overline{k}: \mathsf{NNum}^\dagger(k)_F \to \mathsf{NNum}^\dagger(\overline{k})_F$; the absolute Galois group ${\rm Gal}(\bar k/k)$ is obtained from the Tannakian formalism applied to the category of noncommutative Artin motives $\mathsf{NAM}(k)_F$; and finally $I$ is induced by the inclusion of the latter category in $\mathsf{NNum}^\dagger(k)_F$.

\smallskip

Deligne-Milne's proof \cite{DelMil} of the classical commutative counterpart of \eqref{eq:base-change1} makes full use of ``commutative arguments'' which don't seem to admit noncommutatice analogues. The proof of \eqref{eq:base-change1} is not only very different but moreover much more conceptual from a categorical viewpoint. By extracting the key ingredients of this latter proof we have established in \cite[Appendix A]{MaTa6} a general result about short exact sequences of Galois groups. This led to a new proof of Deligne-Milne's short exact exact sequence which circumvents their ``commutative arguments''.

\medskip
\section{From noncommutative motives to motives via Jacobians}
We have described in \S\ref{sub:bridge} a bridge from motives to noncommutative motives. One can ask is there a bridge in the opposite direction, associating a ``commutative 
shadow" to every noncommutative motive?
This (vague) idea can be implemented using the theory of Jacobians,
suitably extended to the noncommutative world. Let $k \subseteq \bC$ be an algebraically closed field.

\smallskip
\subsection{Jacobians}
Jacobians $J(C)$ of curves $C$ are geometric models for the cohomology $H^1(C)$. The study of Jacobians is in fact one of the historic precursors of the theory of motives. 
Given an arbitrary smooth projective $k$-scheme $X$, the Picard $\mathrm{Pic}^0(X)$ and the Albanese $\mathrm{Alb}(X)$
varieties provide, in a similar way, geometric models 
for the pieces $H^1(X)$ and $H^{2\mathrm{dim}(X)-1}(X)$.
In what concerns the remaining pieces of the cohomology, Griffiths' intermediate Jacobians 
\begin{eqnarray*}
J_i(X) := \frac{H_B^{2i+1}(X,\C)}{F^{i+1}H_B^{2i+1}(X,\C) + H_B^{2i+1}(X,\Z)} && 0 < i < \mathrm{dim}(X)
\end{eqnarray*}
where $H_B$ stands for Betti cohomology and $F$ for the Hodge filtration, are not algebraic. Nevertheless, they contain an algebraic part $J_i^a(X) \subseteq J_i(X)$ defined by the image of the
Abel-Jacobi map $ AJ_i : CH^{i+1}_\Z(X)^{\mathrm{alg}} \to J_i(X)$, where $CH^{i+1}_\Z(X)^{\mathrm{alg}}$ is the group of algebraically trivial cycles of
codimension $i+1$.

\smallskip
\subsection{Pairings}
Given a smooth projective $k$-scheme $X$, consider the following $k$-vector spaces
\begin{eqnarray}\label{eq:NH}
NH_{dR}^{2i+1}(X):= \sum_{C,\gamma_i} \mathrm{Image} \big(H^1_{dR}(C) \stackrel{H^1_{dR}(\gamma_i)}{\longrightarrow} H^{2i+1}_{dR}(X) \big) && 0 \leq i \leq \mathrm{dim}(X)-1\,,
\end{eqnarray}
where $C$ is a smooth projective curve and $\gamma_i: h_\bQ(C) \to h_\bQ(X)(i)$ is a morphism in $\mathsf{Chow}_\bQ(k)$. Intuitively speaking, \eqref{eq:NH} are the odd pieces of de Rham cohomology that are generated by curves. By restricting the classical intersection pairings on de Rham cohomology to \eqref{eq:NH} one obtains
\begin{eqnarray}\label{eq:pairings1}
\langle-,- \rangle : NH_{dR}^{2\mathrm{dim}(X)-2i-1}(X) \times NH_{dR}^{2i+1}(X) \longrightarrow k && 0 \leq i \leq \mathrm{dim}(X)-1\,.
\end{eqnarray}
\subsection{Jacobians of noncommutative motives}
In \cite[Theorem 1.3]{MaTa5} the authors constructed a $\bQ$-linear additive ``Jacobian'' functor
\begin{equation}\label{eq:Jacobian}
{\bf J}(-): \mathsf{NChow}(k)_\bQ \longrightarrow \mathrm{Ab}_\bQ(k)
\end{equation}
with values in the category of abelian varieties up to isogeny. Among other properties, one 
has an isomorphism ${\bf J}(\mathrm{perf}_{\mathsf{dg}}(X)) \simeq \prod_{i=0}^{\mathrm{dim}(X)-1} J_i^a(X)$ whenever the above pairings \eqref{eq:pairings1} 
are non-degenerate. As explained in loc. cit., this is always the case for $i=0$ and $i=\mathrm{dim}(X)-1$ and the remaining cases 
follow from Grothendieck's standard conjecture of type $B$. Hence, the pairings \eqref{eq:pairings1} are non-degenerate 
for curves, surfaces, abelian varieties, complete intersections, uniruled threefolds, rationally connected fourfolds, 
and for any smooth hypersurface section, product, or finite quotient thereof (and if one trusts Grothendieck for all smooth projective $k$-schemes).

\smallskip 

Given a noncommutative Chow motive $N$, the abelian variety ${\bf J}(N)$ was constructed as follows: 
\begin{itemize}
\item first, via $h^1_\bQ(-)$ and the fully-faithful functor $R_{\mathsf{N}}$ (see \S\ref{sub:NCnumerical}), one observes that $\mathrm{Ab}_\bQ(k)$ identifies with a semi-simple abelian full subcategory of $\mathsf{NNum}(k)_\bQ$;
\item secondly, the semi-simplicity of $\mathsf{NNum}(k)_\bQ$ implies that $N$ admits a unique direct sum decomposition $S_1 \oplus \cdots \oplus S_n$ into simple objects; 
\item finally, one defines ${\bf J}(N)$ as the smallest piece of the noncommutative numerical motive $N \simeq S_1 \oplus \cdots \oplus S_n$ which contains the simple objects belonging to $\mathrm{Ab}_\bQ(k)$.
\end{itemize}
Roughly speaking, \eqref{eq:Jacobian} show us that the classical theory of Jacobians can be extended to the noncommutative world as long as one works with all the intermediate Jacobians simultaneously. Note that this restriction is an intrinsic feature of the noncommutative world which cannot be avoided because as soon as one passes from $X$ to $\mathrm{perf}_{\mathsf{dg}}(X)$ one loses track of the individual pieces of $H^\ast_{dR}(X)$ (see \eqref{eq:HKR}).

\subsection{Some applications}
The above theory of Jacobians of noncommutative motives allowed categorical Torelli theorems, a new proof of a classical theorem of Clemens-Griffiths concerning blow-ups of threefolds, and several new results on quadric fibrations and intersections of quadrics; see \cite{MarTab2}. Recently, this theory allowed also the proof of a conjecture of Paranjape \cite{Paranjape} in the case of a complete intersection of either two quadrics or three odd-dimensional quadrics; see \cite{MarTab1}.

\section{Applications to motivic decompositions and to Dubrovin's conjecture}\label{ExCollSec}

\smallskip
\subsection{Motivic decompositions}
It is well-know that $ h(\P^n)= 1 \oplus \bL \oplus \cdots \oplus \bL^{\otimes n}$; see Manin \cite{Man}. Other examples of motivic decompositions containing only $\otimes$-powers of the Lefschetz motive arise from quadrics.
Given a non-degenerate quadratic form $(V,q)$ of dimension $n\geq 3$ defined over an algebraically closed field $k$ of characteristic zero, let $Q_q \subset \bP(V)$
be the associated smooth projective quadric of dimension $d:=n-2$. The 
motivic decomposition of $Q_q$ (due to Rost \cite{Rost}) is the following 
$$
h_\Q(Q_q) \simeq \left\{ \begin{array}{lcl}
1 \oplus {\bL} \oplus \cdots \oplus {\bL}^{\otimes d}  &&  d\,\,\text{odd} \\
1 \oplus {\bL} \oplus  \cdots \oplus {\bL}^{\otimes d} \oplus {\bL}^{\otimes(d/2)}  && 
d \,\,\text{even} \,.
\end{array} \right.
$$
Fano 3-folds also fit in this pattern. In this case, thanks to Gorchinskiy-Guletskii \cite{GoGu}, we have
$$ h_\Q(X) \simeq 1 \oplus h_\Q^1(X) \oplus {\bL}^{\oplus b} \oplus (h_\Q^1(J)\otimes {\bL}) \oplus ({\bL}^{\otimes 2})^{\oplus b} \oplus h_\Q^5(X) \oplus {\bL}^{\otimes 3}\,,$$
where $h_\Q^1(X)$ and $h_\Q^5(X)$ are the Picard and Albanese motives, $b=b_2(X)=b_4(X)$, and $J$ is an abelian variety which is isogenous to the intermediate Jacobian when $k=\C$. Whenenver the odd cohomology of $X$ vanishes, this motivic decomposition
reduces to a direct sum of $\otimes$-powers of the Lefschetz motive. Further exemples include toric varieties and certain homogeneous spaces (see Brosnan \cite{Bro}), and moduli spaces of pointed curves of genus zero (see Chen-Gibney-Karshen \cite{CGK}).

\smallskip
\subsection{Full exceptional collections}\label{sub:exceptional}

A collection of objects $\{ \cE_1,\ldots, \cE_n \}$ in a $k$-linear triangulated category $\cC$ is called  exceptional if ${\rm RHom}(\cE_i,\cE_i)=k$, for all $i$, and ${\rm RHom}(\cE_i,\cE_j)=0$ for all $i>j$. 
It is called full if the objects $\cE_1, \cdots, \cE_n$ generate the triangulated category $\cC$.

\smallskip

The derived category $\mathrm{Perf}(X)\simeq\cD^b(\mathrm{Coh}(X))$ of a smooth projective $k$-scheme $X$ admits a full exceptional collection in several cases: projective spaces (see Beilinson \cite{Beil}), quadrics (see Kapranov \cite{Kapr}), toric varieties (see Kawamata \cite{Kawa}), certain homogeneous spaces 
(see Kuznetsov-Polishchuk \cite{KuPo}), moduli spaces of 
pointed curves of genus zero (see Manin-Smirnov \cite{ManSmi}), 
and Fano 3-folds with vanishing odd cohomology (see Ciolli \cite{Ciolli}). In all these examples the corresponding motivic decomposition contain only $\otimes$-powers of the Lefschetz motive. It is therefore natural to ask if there is a relation between these two notions. The answer is ``yes'', as we now explain.

\smallskip
\subsection{From exceptional collections to motivic decompositions}
Let $X$ be a smooth projective $k$-scheme for which $\mathrm{perf}(X)$ admits a full exceptional collection $\{ \cE_1,\ldots, \cE_n \}$. Then, as proved in \cite[Theorem 1.3]{MaTa4}, there is a choice of integers $\ell_1,\ldots, \ell_n \in \{0, \ldots, \mathrm{dim}(X)\}$ such that 
\begin{equation}\label{eq:motivic-decomposition}
h_\Q(X) \simeq \bL^{\otimes \ell_1} \oplus \cdots \oplus \bL^{\otimes \ell_n}\,.
\end{equation}
The motivic decomposition \eqref{eq:motivic-decomposition} was obtained as follows: 
\begin{itemize}
\item first, as mentioned in \S\ref{sub:computations}, the noncommutative Chow motive $U(\mathrm{perf}_{\mathsf{dg}}(X))_\bQ$ decomposes into the direct sum of $n$ copies of  $U(k)_\bQ$; 
\item secondly, the commutativity of diagram \eqref{eq:diag-main} implies that $h_\bQ(X)$ (considered as an object of the orbit category  $\mathsf{NChow}_\Q(k)/_{\!\!-\otimes \Q(1)}$) decomposes into the direct sum of $n$ copies of ${\bf L}$; 
\item finally, one observes that the ``fiber'' of ${\bf L}$ under the projection functor from noncommutative Chow motives to its orbit category consists solely of $\otimes$-powers of the Lefschetz motive.
\end{itemize}
\smallskip

The decomposition \eqref{eq:motivic-decomposition} holds more generally with $X$ a smooth proper Deligne-Mumford stack.

\smallskip

The decomposition \eqref{eq:motivic-decomposition} has recently greatly refined: instead of working with $\bQ$-coefficients it suffices to invert the prime factors of the integer $(2\mathrm{dim}(X))!$; consult \cite{MarTab3} for details.

\smallskip

\subsection{Dubrovin's conjecture}
At his ICM address \cite{Dubrovin}, Dubrovin conjectured a striking connection between Gromov-Witten invariants and derived categories of coherent sheaves. The most recent formulation, due to Hertling-Manin-Teleman \cite{HMT}, is the following:

\smallskip
 
{\it Conjecture: Given a smooth projective $\bC$-scheme $X$, the following holds: 
\begin{itemize}
\item[(i)] The quantum cohomology of $X$ is semi-simple if and only if $X$ is Hodge-Tate (i.e the Hodge numbers $h^{p,q}(X)$ are zero for $p \neq q$) and $\mathrm{perf}(X)$ admits a full exceptional collection;
\item[(ii)] The Stokes matrix of the structure connection of the quantum cohomology identifies with the Gram matrix of the full exceptional collection.
\end{itemize}
}
 
\smallskip 
 
Thanks to the work of Bayer, Golyshev, Guzzeti, Ueda, and others (see \cite{Bayer,Golyshev,Guzzetti,Ueda}), items (i)-(ii) are known to be true in the case of projective spaces (and its blow-ups) and Grassmannians. Item (i) also holds for minimal Fano threefolds. Moreover, Hertling-Manin-Teleman proved in \cite{HMT} that the Hodge-Tate property follows from the semi-simplicity of quantum cohomology. 

\smallskip 
 
Making use of the above motivic motivic decomposition \eqref{eq:motivic-decomposition}, the authors proved in \cite[Proposition 1.9]{MaTa4} that the Hodge-Tate property follows also from the existence of a full exceptional collection. As a consequence, this assumption can be removed from item (i) of Dubrovin's conjecture.

\end{document}